\documentclass{article}
\usepackage[margin=1in]{geometry}

\usepackage{graphicx} 
\usepackage{amsmath,amsthm,amssymb}
\usepackage{natbib}
\newtheorem{proposition}{Proposition}[section]

\newtheorem{corollary}{Corollary}[section]
\newtheorem{definition}{Definition}[section]

\newtheorem{example}{Example}[section]
\newtheorem{remark}{Remark}[section]

\title{Dependence properties between marginals of discrete-time maximum processes}
\author{J. Baz and F. Pellerey}
\date{}

\begin{document}

\maketitle

\begin{abstract}
Consider a discrete time sequence $(X_n,n\in\mathbb{N})$ and the corresponding process $(M_n,n\in\mathbb{N})$ which, at any time $n\in\mathbb{N}$, takes the maximum value reached by the process $(X_n,n\in\mathbb{N})$ up to that time. Clearly, there exists a dependence  between the values assumed by the process $(M_n,n\in\mathbb{N})$ at two distinct points in time. This paper formally describes such a dependence in terms of the copula of the two random values. Various scenarios with respect to the sequence $(X_n,n\in\mathbb{N})$ are considered, e.g., without assuming independence in its components and/or identical distribution. Such explicit expressions of the copula are then used to provide results that describe how the dependence between the maxima reached by the sequence $(X_n,n\in\mathbb{N})$ in different times varies when both the time interval between them and/or their positions in the sequence are increased. In this context, asymptotic results are also provided.    
\end{abstract}

\section{Introduction}

Let $(X_n,n\in\mathbb{N})$ be sequence of random variables indexed in discrete time. In many scenarios such as climate studies~\citep{wergen2010record,jouni2025exploring}, sports~\citep{gembris2002trends,tang2024analysis} or economy~\citep{wergen2011record,kim2021statistical}, keeping track of the current value of the maximum of the sequence is of interest. In this regard, let us define the sequence $(M_n,n\in\mathbb{N})$ letting
\begin{equation*}
    M_n=\max(X_1,\dots,X_n),
\end{equation*}
i.e., the maximum between the first $n$ values of the original sequence. We devote this paper to the study of dependence properties of this type of sequences. Under different assumptions for the initial sequence $(X_n,n\in\mathbb{N})$, we will compute the copula of pairs of variables in $(M_n,n\in\mathbb{N})$, proving positive dependence properties, deriving dependence stochastic inequalities and studying the asymptotic behavior.

It is worth to point out the connection of the sequence $(M_n,n\in\mathbb{N})$ with similar processes often considered in the literature, such as in record theory~\citep{nevzorov1988records}. In that case, the sequence of upper record times $L_1,L_2,\dots$ is defined as $L_1=1$ and $L_{n+1}=\min \{j\mid j\geq L_n, X_j>X_{L_n}\}$ and then the sequence of records is defined as $Y_n=X_{L_n}$. The main difference between $(M_n,n\in\mathbb{N})$ and $(Y_n,n\in\mathbb{N})$ is that, while $M_n$ saves the value of the current maximum in each time (so repeated values appear), $Y_n$ only adds a new value when an update of the maximum is made. For instance, given the sequence $(1,2,2,1,4,3)$, the maximum process is given by $(1,2,2,2,4,4)$, while the associated record values are just $(1,2,4)$. One can find some few papers dealing with dependence inequalities between records of some  models~\citep{khaledi2005dependence,wergen2011correlations}, some of them as particular cases of Generalized Ordered Statistics~\citep{kamps1995concept}. 

The maximum process is also related to the study of order statistics with a variable sample size. Given a sequence of random variables $(X_n,n\in\mathbb{N})$, we will denote the $i$-th order statistic associated with $(X_1,\dots,X_n)$, i.e. the $i$-th smallest value among the first $n$ variables, as $X_{(i:n)}$. One can find some dependence comparison results in the literature (see~\cite{averous2005dependence,navarro2010relationships,kochar2012stochastic}). However, the compared vectors of order statistics are obtained for fixed, maybe different between them, sample sizes. In addition, only the case of independent and identically distributed random variables is commonly considered.

Finally, the study presented here could be also connected with some papers dealing with dependence ordering of Markov chains~\citep{fang1994decrease,hu2000comparisons,kulik2004monotonicity}. It is easy to see that, since one can write $M_{n+1}=\max(M_n,X_{n+1})$, the maximum process is Markov if the initial random variables are independent. However, most of such results impose the stationarity of the sequence, while the maximum process is almost surely increasing, or impose a finite number of states for the considered Markov chain. 

This paper, although connected with the previously mentioned works, deals with a more general framework in which a bigger family of initial sequences $(X_n,n\in\mathbb{N})$ are considered, including exchangeable sequences or independent but differently distributed random variables. Given two different times $n_1,n_2\in\mathbb{N}$, we derive the copula of $(M_{n_1},M_{n_2})$ and, considering two additional times  $n_3,n_4\in\mathbb{N}$, we give conditions for the inequality $(M_{n_1},M_{n_2})\preceq_c (M_{n_3},M_{n_4})$, where $\preceq_c$ stands for the concordance dependence order. This dependence order compares the dependence structure of the involved random vectors, so our results allow to identify pairs of random variables in $(M_n,n\in\mathbb{N})$ with a stronger positive dependence.

The results given in this paper, beyond their own self-interest, are relevant in different scenarios. Firstly, we can find in the literature some papers studying the copula of the maximum and of pairs of random variables~\citep{li2011generalized}, or the minimum and maximum of random variables under some particular conditions~\citep{vidmar2015copulas}. Here, we study copulas considering the maximum of any number of variables and under simple and natural conditions. Secondly, we provide rigorous proofs to intuitive facts regarding stochastic inequalities and asymptotic behavior of the dependence in the maximum process. For instance, one expects a stronger dependence when considering the maximum process in two closer time points and, when fixing the distance in time, a stronger dependence as the time increases. Finally, it has to be pointed out that the convergence of $(M_n,n\in\mathbb{N})$, under the correct transformation, has been studied extensively in the so-called extreme value theory~\citep{leadbetter2012extremes, ji2023extreme}. The asymptotic results for copulas given in this paper could be used to extend these studies to the convergence of bivariate distributions in the maximum process rather than just univariate distributions. These kind of results are for example provided in Section 3 of~\cite{leadbetter2012extremes}, but making use of a measure of ``distance from independenc'' between the values assumed by the maximum process at different times. The results presented here more clearly describe such properties of monotonicity in dependence, being based on the copula of the two values. Another potential application involves the computation of some conditional probabilities for the case of Archimedean copulas, see Section~\ref{sec:archimedean}, since this type of copulas have been appeared in recent papers dealing with extreme values~\citep{das2024ordering,shrahili2024stochastic}.

Throughout the paper, the marginal distribution functions of the random variables are assumed to be continuous. Strict monotonicity is not required; inverses of distribution functions are understood as generalized inverses, as specified in the next section.

\section{Preliminaries}

The main tool we are going to use to study the dependence structure of maximum processes are copulas. Copulas are multivariate distributions with standard uniform marginals and are widely used to model the dependence structure of random vectors~\citep{nelsen2006introduction}. In particular, the distribution of any random vector can be decomposed in the marginal distributions, how the variables behave separately, and its copula, how they behave together.

\begin{proposition}\citep{durante2016principles,sklar1959fonctions}\label{thm:sklar}
Let $\mathbf{X}$ be a random vector with marginal distribution functions $F_1,\dots,F_n$ and joint distribution function $F$. Then, there exists a copula $C:[0,1]^n\to[0,1]$ such that
\begin{equation*}
F(x_1,\dots,x_n)=C\left(F_1(x_1),\dots,F_n(x_n)\right),\ \forall\mathbf{x}\in\mathbb{R}^n.
\end{equation*}
In addition, if $F_1,\dots,F_n$ are continuous, the copula is unique.
\end{proposition}

Throughout the paper, all marginal distribution functions are assumed to be continuous. For any distribution function $G$, we use the generalized inverse
\[
G^{-1}(u)=\inf\{x\in\mathbb{R}:G(x)\ge u\}, \qquad u\in(0,1).
\]
If $G$ is continuous, then
\[
G\bigl(G^{-1}(u)\bigr)=u, \qquad u\in(0,1),
\]
even when $G$ is not strictly increasing.

Consequently, if $H$ is a bivariate distribution function with continuous marginal distribution functions $G_1$ and $G_2$, its unique copula satisfies
\[
C(u,v)=H\bigl(G_1^{-1}(u),G_2^{-1}(v)\bigr), \qquad u,v\in(0,1).
\]
Copula formulas involving generalized inverses are initially
understood on $(0,1)^2$. Their boundary values are obtained by continuous extension and satisfy
\[
C(u,0)=C(0,v)=0,\qquad C(u,1)=u,\qquad C(1,v)=v.
\]

Any bivariate copula is bounded by $\max(0,u+v-1)\leq C(u,v)\leq \min(u,v)$. Such bounds are known as the Fréchet–-Hoeffding bounds~\citep{nelsen2006introduction}. Similarly, the multivariate survival function can be decomposed in the marginal survival functions and the so called survival copula~\citep{nelsen2006introduction}. In this work, the case of Archimedean copulas will be of special interest. They are defined in terms of a generator $\phi$, that is used to construct the copula for the adequate dimension. We recall that a completely monotonic function $f:D\to\mathbb{R}$ is a real continuous function such that their derivatives exists and are of alternating sign, i.e. $(-1)^k(\text{d}^k/\text{d}x^k)f(x)\geq 0$ for any $x\in D$ and $k\in\mathbb{N}\cup\{0\}$~\citep{nelsen2006introduction}.

\begin{definition}\citep{nelsen2006introduction}
A copula $C:[0,1]^n\to [0,1]$ is said to be Archimedean if there exists a continuous, strictly decreasing function $\phi:[0,1]\to [0,\infty]$ fulfilling $\phi(1)=0$, with completely monotonic pseudo-inverse $\phi^{-1}$ and such that
\begin{equation*}
C(\mathbf{u})=\phi^{-1}\left(\phi(u_1)+\dots+\phi(u_n)\right).
\end{equation*}
\end{definition}

As particular case, we will consider Clayton copulas, for which $\phi(x)=(x^{-\theta}-1)/\theta$ with $\theta\geq 0$. In addition, the Marshall-Olkin family, which is not Archimedean, will appear naturally when working with independent and identically distributed random variables.

\begin{definition}\label{def:MarshallOlkin}\citep{nelsen2006introduction}
A bivariate copula $C$ is said to be a Marshall-Olkin copula if it can be expressed as
\begin{equation*}
    C(u,v)=\min(u^{1-\alpha}v,uv^{1-\beta})
\end{equation*}
with $\alpha,\beta\in [0,1]$.
\end{definition}

In general, if the copula takes greater values, there is a stronger positive dependence between the associated random variables. For the comparison of the dependence structure of two random vectors, we can consider the pointwise comparison of copulas. In this direction, the concordance order will be used. Usually, this order is defined between copulas, although it can be extended trivially for random vectors with continuous components since they have a unique copula.

\begin{definition}\citep{nelsen2006introduction}
Let $\mathbf{X}$ and $\mathbf{Y}$ be two bivariate random vectors with continuous components and copulas, respectively, $C_x$ and $C_y$. Then, it is said that $\mathbf{X}$ is smaller than or equal to $\mathbf{Y}$ in the concordance order if $C_x(\mathbf{u})\leq C_y(\mathbf{u})$ for any $\mathbf{u}\in [0,1]^2$. It will be denoted as $\mathbf{X}\preceq_c\mathbf{Y}$.
\end{definition}

If the involved random vectors have more than two components, the concordance order also requires the ordering of the survival copulas, see~\cite{nelsen2006introduction}. The concordance order implies the inequality of well-known dependence measures such as Kendall’s tau, Spearman’s rho and Blomqvist’s $\beta$~\citep{nelsen2006introduction}.

Throughout the paper, we will occasionally refer to the concept of (usual) stochastic order. For this reason, we provide its definition here.

\begin{definition}\citep{shaked2007stochastic}
Let $\mathbf{X}$ and $\mathbf{Y}$ be two random vectors. $\mathbf{X}$ is said to be smaller than or equal to $\mathbf{Y}$ in the usual stochastic order, denoted as $\mathbf{X}\leq_{st}\mathbf{Y}$, if $P(\mathbf{X}\in U)\leq P(\mathbf{Y}\in U)$ for any measurable upper set $U\subseteq \mathbb{R}^n$.
\end{definition}

It is easy to observe that, in the univariate case, the usual stochastic order is characterized by the pointwise comparison of the distribution functions, i.e. $X\leq_{st}Y\iff F_X(t)\geq F_Y(t)\ \forall t\in \mathbb{R}$~\citep{shaked2007stochastic}.

\section{The maximum process}

Let $(X_n,n\in\mathbb{N})$ be a sequence of random variables. As stated in the introduction, its associated maximum process is defined as $(M_n,n\in\mathbb{N})$ where $M_n=\max(X_1,\dots,X_n)$. Trivially, the maximum process is almost surely increasing, i.e. $P(M_{n+1}\geq  M_{n})=1$. It is also interesting to note that, if $(X_n,n\in\mathbb{N})$ consists on independent random variables, then $M_{n+1}$ is conditionally independent of $(M_1,\dots,M_{n-1})$ given the value of $M_n$ (since $M_{n+1}=\max(M_n,X_{n+1})$ and $X_{n+1}$ is independent of $(X_1,\dots,X_{n-1})$, thus also of $(M_1,\dots,M_{n-1})$). That is, the maximum process is a Markov chain. Since each $X_i$ has a continuous distribution, for every $t\in\mathbb{R}$,
\[
P(M_n=t)\le\sum_{i=1}^{n}P(X_i=t)=0.
\]

Thus, each $M_n$ has a continuous distribution function, although this function need not be strictly increasing. In particular, every pair $(M_n,M_{n+k})$ has a unique copula.

There are two different sources of dependence between the variables in the maximum process. First, we have a natural positive dependence given by the definition of the sequence. Since the maximum is an increasing function, we expect large values of the maximum process at any given times should imply large values at subsequent times. In addition, we also expect variables that are closer in time to have a stronger dependence. The second effect is the initial dependence structure of the sequence $(X_n,n\in\mathbb{N})$. In general, a stronger positive dependence of $(X_n,n\in\mathbb{N})$ should imply a stronger positive dependence of $(M_n,n\in\mathbb{N})$. We will focus mainly in the dependence due to the first source, the construction of the maximum process, while the other source will be relevant only in Proposition~\ref{prop:conc_dependence} and its discussion. Let us start by providing a formula for the copula of the maximum process in a very general case (under the mentioned assumption on the distribution of the $X_i$).

\begin{proposition}\label{prop:inicial}
Let $(X_n,n\in\mathbb{N})$ be a sequence of random variables having continuous cumulative distributions. Denote as $F_{(n:n)}$ the distribution function of $M_n=\max(X_1,\dots,X_n)$ and, for integers $m>n$ and values of $x$ such that $P(M_n\le x)>0$, let $F_{\{n+1:m\},x}$ denote the distribution function of
$
\bigl[\max(X_{n+1},\ldots,X_m)\mid M_n\le x\bigr].
$
Then,
\begin{equation*}
    C_{n,n+k}(u,v)=\min\left(uF_{\{n+1:n+k\}, F^{-1}_{(n:n)}(u)}(F^{-1}_{(n+k:n+k)}(v)),v\right)
\end{equation*}
is the copula of $(M_n,M_{n+k})$.
\end{proposition}

\begin{proof} Let us define $Y=\max(X_{n+1},\ldots,X_{n+k})$, so that $M_{n+k}=\max(M_n,Y)$, and then define
\[H(x,y)=P({M_n}\le x,Y\le y)=F_{M_n}(x)\,P(Y\le y\mid {M_n}\le x),
\]
whenever $F_{M_n}(x)>0$. The joint distribution function of $(M_n,M_{n+k})$ is therefore $G(x,y)=P({M_n}\le x,\max({M_n},Y)\le y)$. We distinguish two cases:
\begin{itemize}
\item If $x\le y$, then $G(x,y)=H(x,y)$. Moreover,
\[
H(x,y)\le P({M_n}\le y,Y\le y)=F_{M_{n+k}}(y).
\]

\item If $x>y$, then $G(x,y)=F_{M_{n+k}}(y)$.
By inclusion of events,
\[
H(x,y)\ge P({M_n}\le y,Y\le y)=F_{M_{n+k}}(y).
\]
\end{itemize}

In both cases we have
\[
G(x,y)=\min\left\{F_{M_n}(x)\,P(Y\le y\mid {M_n}\le x),\,F_{M_{n+k}}(y)\right\}.
\]

The marginal distribution functions $F_{M_n}$ and $F_{M_{n+k}}$ are continuous. For $u,v\in(0,1)$, set
\[
x=F_{M_n}^{-1}(u),\qquad y=F_{M_{n+k}}^{-1}(v),
\]
where the inverses are generalized inverses. Continuity gives
\[
F_{M_n}(x)=u,
\qquad
F_{M_{n+k}}(y)=v.
\]
In particular, $F_{M_n}(x)=u>0$, so the conditional distribution
appearing in the formula is well defined. By Sklar's theorem,
substituting these values into the preceding expression for
$G(x,y)$ yields
\[
C_{n,n+k}(u,v)
=
\min\left\{
u\,F_{\{n+1:n+k\},F_{M_n}^{-1}(u)}
\left(F_{M_{n+k}}^{-1}(v)\right),
\,v
\right\}.
\]
The substitution applies directly for $u,v\in(0,1)$, while the boundary values are obtained by continuous extension.
\end{proof}

While valid, such a formula is not really operative since it involves a conditional distribution that is not easy to compute. An alternative formula for the copula of the maximum process can be found in Remark~4 in~\cite{vidmar2015copulas}. However, the formula is not operative either and it holds under more strict requirements, see Assumption~1 in~\cite{vidmar2015copulas}. In the following sections, we reach simpler formulas for the copulas by restricting the initial random sequence $(X_n,n\in\mathbb{N})$ to have identically distributed or independent variables.

\subsection{Maximum processes of independent random variables}

Consider $(X_n,n\in\mathbb{N})$ to be a sequence of independent random variables, with possibly different distributions.
Independence allows us to simplify the formula of Proposition~\ref{prop:inicial}.

\begin{proposition}\label{thm:copula_ind}
Let $(X_n,n\in\mathbb{N})$ be a sequence of independent random variables
with continuous distribution functions, not necessarily
identically distributed. Denote the distribution function of $M_n$ as $G_n$, with $n\in\mathbb{N}$. Then,
\begin{equation*}
    C_{n,n+k}(u,v)=\min\left(\frac{uv}{G_n(G^{-1}_{n+k}(v))},v\right)
\end{equation*}
is the copula of $(M_n,M_{n+k})$.
\end{proposition}

\begin{proof}
Let $F_i$ denote the distribution function of $X_i$. By
independence,
\[
G_m(x)=\prod_{i=1}^{m}F_i(x),
\]
so each $G_m$ is continuous. Fix now $u,v\in(0,1)$ and let
\[
a=G_n^{-1}(u),
\qquad
b=G_{n+k}^{-1}(v).
\]
By continuity, $G_n(a)=u$ and $G_{n+k}(b)=v$. Moreover, since $M_n\le M_{n+k}$ almost surely,
\[
G_n(b)\ge G_{n+k}(b)=v>0.
\]

Using now independence, we obtain
\begin{align*}
C_{n,n+k}(u,v)
&=P(M_n\le a,M_{n+k}\le b)=
\prod_{i=1}^{n}F_i(\min\{a,b\})
\prod_{j=n+1}^{n+k}F_j(b)=
G_n(\min\{a,b\})\,
\frac{G_{n+k}(b)}{G_n(b)}\\
&=
\min\{u,G_n(b)\}\,\frac{v}{G_n(b)}=
\min\left\{\frac{uv}{G_n(b)},v\right\}.
\end{align*}
Substituting $b=G_{n+k}^{-1}(v)$ one gets the formula on $(0,1)^2$. The boundary values follow by continuity.
\end{proof}

The latter copula is similar to the copulas considered in~\cite{li2011generalized}. In particular, if one replaces the distribution functions $G_n$ and $G_{n+k}$ by a pair of survival functions, the resulting function is a survival copula, as in~\cite{li2011generalized}. Notice that the copula only depends on the distribution functions of $M_n$ and $M_{n+k}$. Using this fact, we can immediately prove the following stochastic inequality.

\begin{proposition}\label{prop:conc_ind}
Let $(X_n,n\in\mathbb{N})$ be a sequence of independent random variables having continuous distributions, non necessarily identically distributed. Then, $(M_n,M_{n+k})\succeq_c\allowbreak(M_n,M_{n+k'})$ for any $n,k,k'\in\mathbb{N}$ such that $k\leq k'$.
\end{proposition}

\begin{proof}
Both copulas are
\begin{equation*}
    C_{n,n+k}(u,v)=\min\left(\frac{uv}{G_n(G^{-1}_{n+k}(v))},v\right),\quad C_{n,n+k'}(u,v)=\min\left(\frac{uv}{G_n(G^{-1}_{n+k'}(v))},v\right),
\end{equation*}
being the only difference the quantile functions $G^{-1}_{n+k'}$ and $G^{-1}_{n+k}$. Since $M_{n+k}\leq_{a.s.}M_{n+k'}$, thus also $M_{n+k}\leq_{st}M_{n+k'}$, $G^{-1}_{n+k}(p)\leq G^{-1}_{n+k'}(p)$ for any $p\in(0,1)$. Therefore, since $G_n$ is increasing, the result holds.
\end{proof}

\begin{figure}
    \centering
    \includegraphics[width=0.6\linewidth]{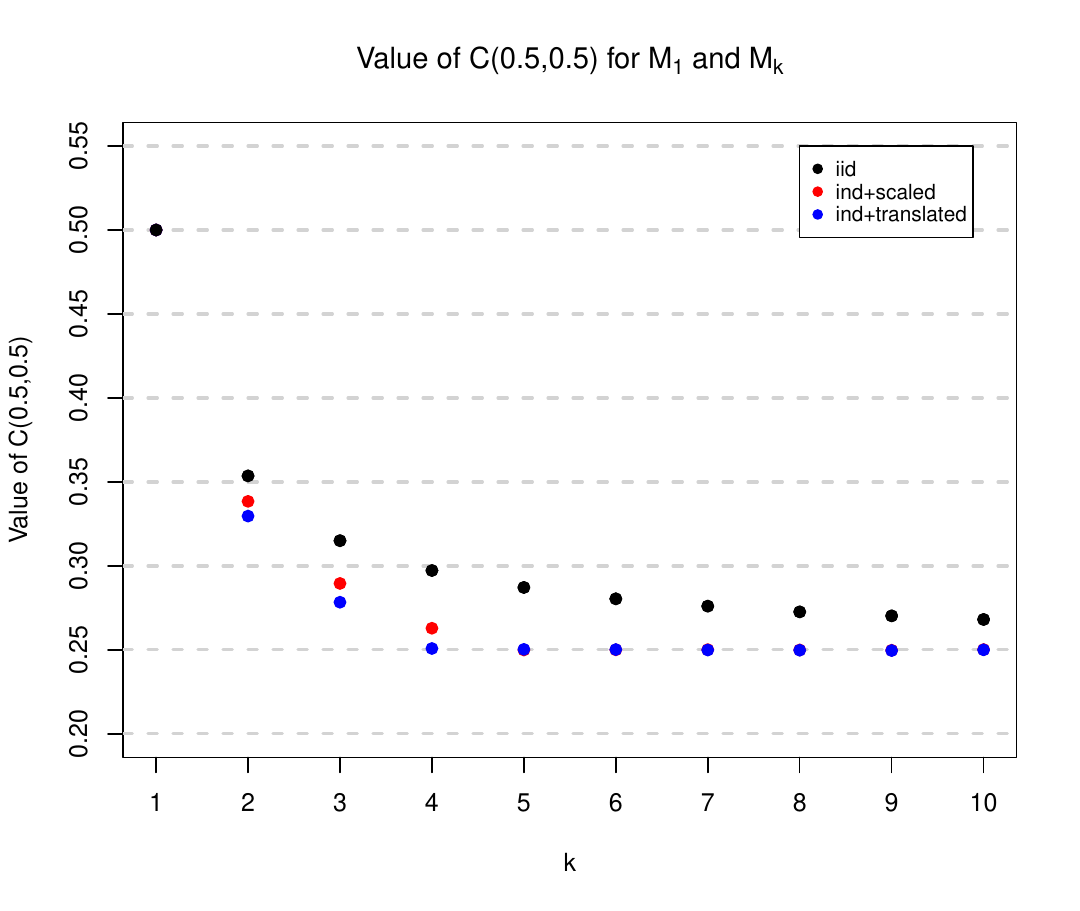}
    \caption{Simulated values of the copula of $(M_1,M_{k})$ in $(0.5,0.5)$ for a sequence of independent random variables $(X_n,n\in\mathbb{N})$ such that $X_n$ is uniform over $[0,1]$ (black), $X_n$ has uniform distribution on $[0,1+n/10]$ (red) and $X_n$ has uniform distribution on $[n/10,1+n/10]$ (blue), using a sample size of $10^6$. Points for the second and third cases are overlapped for $k\geq 5$.}
    \label{fig:1}
\end{figure}

Essentially, the latter results says that the dependence is stronger if the components are closer in time. An illustration of such monotonicity can be found in Figure~\ref{fig:1}.

If one assumes $(X_n,n\in\mathbb{N})$ having not only independent variables, but also identically distributed, then the expression for the copulas becomes very simple. In fact, assuming that the random variables have the same distribution, the copula takes a simple expression that does not depend on the marginals.

\begin{corollary}\label{cor:iid_max}
Let $(X_n,n\in\mathbb{N})$ be a sequence of independent and identically distributed random variables having continuous distribution. Then,
\begin{equation*}
    C_{n,n+k}(u,v)=\min\left(uv^{\frac{k}{n+k}},v\right)
\end{equation*}
is the copula of $(M_n,M_{n+k})$.
\end{corollary}

\begin{proof}
Notice that, using independence, we have that  $G_n(G^{-1}_{n+k}(v))=F(F^{-1}(v^{1/(n+k)}))^n=v^{n/(n+k)}$, where $F$ is the common distribution function of the involved variables. Then, use the formula of Proposition~\ref{thm:copula_ind}.
\end{proof}

If  one takes $\alpha=1$ and $\beta=\frac{n}{n+k}$, it follows that
\begin{equation*}
    C_{n,n+k}(u,v)=\min\left(uv^{\frac{k}{n+k}},v\right)=\min\left(u^{1-\alpha}v,uv^{1-\beta}\right),
\end{equation*}
thus the copula between two variables in the maximum process of independent and identically distributed random variables is a Marshall-Olkin copula. With the formula of Corollary~\ref{cor:iid_max}, the following inequality is trivially satisfied.

\begin{corollary}\label{cor:iid_conc}
Let $(X_n,n\in\mathbb{N})$ be a sequence of independent and identically distributed random variables having continuous distribution. Then, $(M_n,M_{n+k})\preceq_c(M_{n'},M_{n'+k})$ for any $n,n',k\in\mathbb{N}$ such that $n\leq n'$.
\end{corollary}

\begin{figure}
    \centering
    \includegraphics[width=0.6\linewidth]{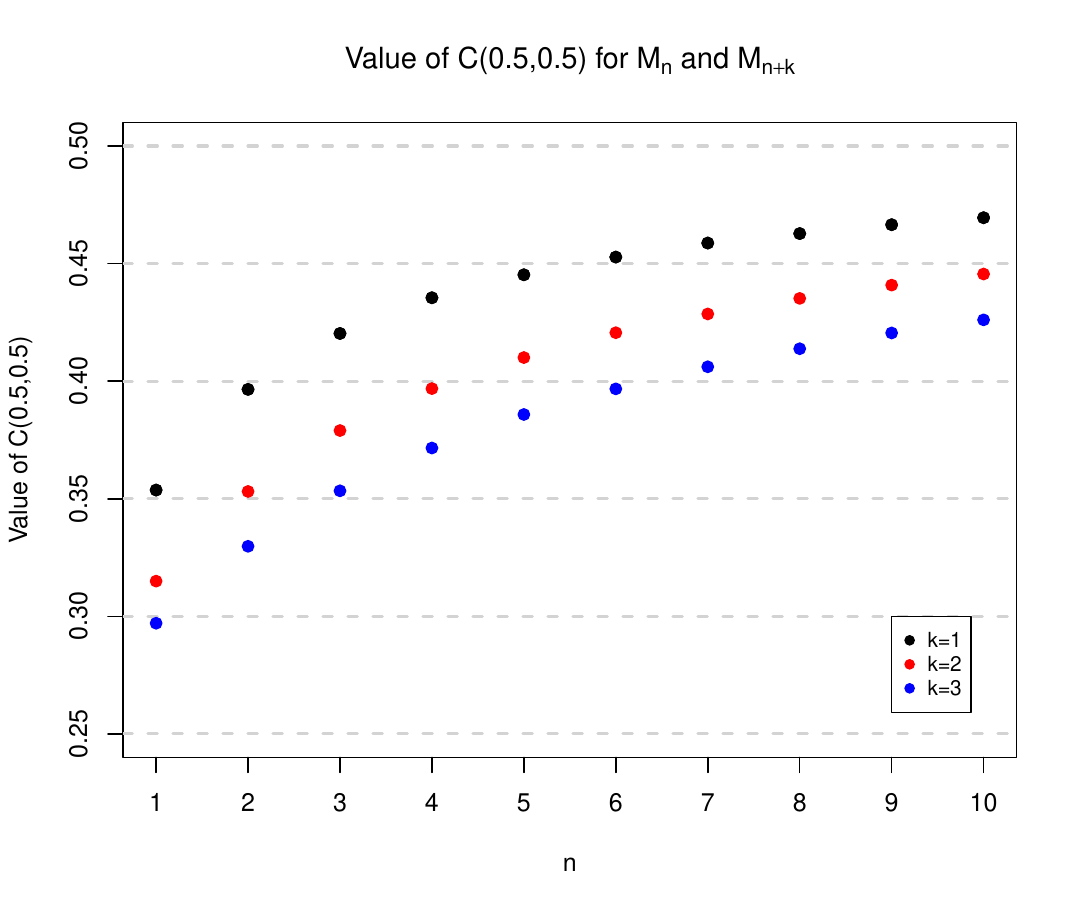}
    \caption{Simulated values of the copula of $(M_n,M_{n+k})$ in $(0.5,0.5)$ for a sequence of independent random variables with uniform distribution on $[0,1]$, using a sample size of $10^6$.}
    \label{fig:2}
\end{figure}

An illustration of the behavior can be found in Figure~\ref{fig:2}

\subsection{Sequences of identically distributed random variables}\label{sec:archimedean}

In order to extend the results provided in the previous section, we need to deal with dependent random variables. Unfortunately, in the general case one can find examples in which the dependence inequality proved in Proposition~\ref{prop:conc_ind} is not longer true.

\begin{example}
Let $X_1$, $X_2$ and $X_3$ be three random variables with uniform distribution over the intervals, respectively, $[0,1]$, $[1,2]$ and $[2,3]$. If $X_1$ and $X_2$ are countermonotone and $X_1$ and $X_3$ are comonotone, i.e. $X_2=_{a.s.}2-X_1$ and $X_3=_{a.s.} 2+X_1$, then we have that $(M_1,M_2)=(X_1,X_2)$ are countermonotone while $(M_1,M_3)=(X_1,X_3)$ are comonotone, thus $(M_1,M_2)\leq_{c}(M_1,M_3)$.
\end{example}

However,  it is possible to obtain some results when the initial sequence has identically distributed marginals but without imposing independence. Let us start computing the formula of the copula in that case.

\begin{proposition}\label{thm:samemarginals_max}
Let $(X_n,n\in\mathbb{N})$ be a sequence of identically distributed random variables having continuous distribution.. Denote the copula of $(X_1,\dots,X_{n+k})$ as $C$ and consider
\begin{equation*}
    \delta_{i}(u)=C(\underbrace{u,\dots,u}_{i},\underbrace{1,\dots,1}_{n+k-i}),
\end{equation*}
with $i\leq n+k$. Also, for $w\in(0,1)$, define $\delta_i^{-1}(w)=\inf\{u\in[0,1]:\delta_i(u)\ge w\}.$ Then,
\begin{equation*}
    C_{n,n+k}(u,v)=\begin{cases}
        v & \text{if } \delta^{-1}_n(u)\geq \delta^{-1}_{n+k}(v),\\
C(\underbrace{\delta^{-1}_n(u),\dots,\delta^{-1}_n(u)}_{n}, \underbrace{\delta^{-1}_{n+k}(v),\dots,\delta^{-1}_{n+k}(v)}_{k}) & \text{if } \delta^{-1}_n(u)< \delta^{-1}_{n+k}(v)
    \end{cases}
\end{equation*}
is the copula of $(M_n,M_{n+k})$.
\end{proposition}

\begin{proof}
Note that each $\delta_i$ is continuous and nondecreasing, with $\delta_i(0)=0$ and $\delta_i(1)=1$, thus $\delta_i\bigl(\delta_i^{-1}(w)\bigr)=w, w\in(0,1).$ Let $F$ denote the common continuous distribution function of the variables $X_i$, and let $G_i$ denote the distribution function of $M_i$. Then $G_i(x)=\delta_i(F(x))$.

Fix $u,v\in(0,1)$ and let $a=\delta_n^{-1}(u)$ and $b=\delta_{n+k}^{-1}(v)$. By continuity of the diagonal sections,
\[
\delta_n(a)=u, \qquad \delta_{n+k}(b)=v.
\]
Moreover, $a,b\in(0,1)$. Define $x=F^{-1}(a)$ and $y=F^{-1}(b)$. Since $F$ is continuous, then $F(x)=a$ and $F(y)=b$ and therefore $G_n(x)=u$ and $G_{n+k}(y)=v$.

By Sklar's theorem and the definition of the maximum process,
\begin{align*}
C_{n,n+k}(u,v)
&=P(M_n\le x,M_{n+k}\le y)=
C(\underbrace{\min\{a,b\},\ldots,\min\{a,b\}}_{n},\underbrace{b,\ldots,b}_{k}),
\end{align*}
where we used $F(\min\{x,y\})=\min\{F(x),F(y)\}=\min\{a,b\}.$

If $a\ge b$, the latter expression equals
\[
C(b,\ldots,b)=\delta_{n+k}(b)=v,
\]
while if $a<b$, it equals
\[
C(\underbrace{a,\ldots,a}_{n},\underbrace{b,\ldots,b}_{k}).
\]
By substituting $a=\delta_n^{-1}(u)$ and $b=\delta_{n+k}^{-1}(v)$, this proves the stated formula on $(0,1)^2$. The boundary
values follow by continuity.
\end{proof}

\begin{remark}\label{remark}
Notice that, under the conditions given in Proposition~\ref{thm:samemarginals_max}, we can write $C_{n,n+k}(u,v)=\min(Q(u,v),v)$, where 
\begin{equation*}
Q(u,v)=C(\underbrace{\delta^{-1}_n(u),\dots,\delta^{-1}_n(u)}_{n}, \underbrace{\delta^{-1}_{n+k}(v),\dots,\delta^{-1}_{n+k}(v)}_{k}).
\end{equation*}
\end{remark}

The fact that $C_{n,n+k}$ is actually a copula can be checked by doing straightforward computations. Notice that, if one considers independence, the formula of Corollary~\ref{cor:iid_max} is obtained. It is also worth to mention that such a copula does not depend on the distribution of the marginals. That is, the copula of the maximum process, if the involved random variables are identically distributed, is independent of what is the election of such a common distribution. This coincides with results involving order statistics, see~\cite{navarro2010relationships}. Latter formula allows us to replicate Proposition~\ref{prop:conc_ind} for sequences considering an Archimedean copula and identically distributed random variables.

\begin{proposition}\label{prop:archi}
Let $(X_n,n\in\mathbb{N})$ be a sequence of independent and identically distributed random variables having continuous distribution such that, for any $n\in\mathbb{N}$,  $(X_1,\dots,X_n)$ has an Archimedean copula with generator $\phi$. Then,
\begin{itemize}
  \item $(M_n,M_{n+k})\succeq_c(M_n,M_{n+k'})$ for any $n,k,k'\in\mathbb{N}$ such that $k\leq k'$;
  \item $(M_{n'},M_{n'+k})\succeq_c(M_n,M_{n+k})$ for any $n,n',k\in\mathbb{N}$ such that $n\leq n'$
\end{itemize}
\end{proposition}

\begin{proof}
Let $C_{n,n+k}$ denote the copula of $(M_n,M_{n+k})$. We first show that
\begin{equation}\label{eq:maximum_archimedean_copula}
C_{n,n+k}(u,v)=
\min\left\{v,\,\phi^{-1}\left(\phi(u)+\frac{k}{n+k}\phi(v)\right)\right\}.
\end{equation}

Fix $u,v\in(0,1)$. Under the Archimedean assumption,
the diagonal sections introduced in Proposition~\ref{thm:samemarginals_max} satisfy
\[\delta_i(t)=\phi^{-1}(i\phi(t)),\qquad\delta_i^{-1}(w)=\phi^{-1}\left(\frac{\phi(w)}{i}\right).
\]
Define now
\[
a=\delta_n^{-1}(u)=\phi^{-1}\left(\frac{\phi(u)}{n}\right),\qquad \textrm{and}\qquad b=\delta_{n+k}^{-1}(v)=\phi^{-1}\left(\frac{\phi(v)}{n+k}\right).
\]
By Proposition~\ref{thm:samemarginals_max} and the Archimedean representation, we have
\[
C_{n,n+k}(u,v)=\phi^{-1}\left(n\phi(\min\{a,b\})+k\phi(b)\right).
\]
Since $\phi$ is decreasing, then it follows
\[
\phi(\min\{a,b\})=\max\left\{\frac{\phi(u)}{n},\,\frac{\phi(v)}{n+k}\right\}.
\]
Consequently,
\begin{align*}
C_{n,n+k}(u,v)
&=
\phi^{-1}\left(n\max\left\{\frac{\phi(u)}{n},\,\frac{\phi(v)}{n+k}\right\}+\frac{k}{n+k}\phi(v)\right)\\&=\phi^{-1}\left(\max\left\{\phi(u)+\frac{k}{n+k}\phi(v),\,\phi(v)\right\}\right)\\&=\min\left\{\phi^{-1}\left(\phi(u)+\frac{k}{n+k}\phi(v)\right),\,v\right\},
\end{align*}
where the last equality follows because $\phi^{-1}$ is also decreasing. This proves \eqref{eq:maximum_archimedean_copula} on $(0,1)^2$; the boundary values follow by continuity. Now let $k\le k'$. Since $k/(n+k)\le k'/(n+k')$ and $\phi(v)\ge0$, we have
\[
\phi(u)+\frac{k}{n+k}\phi(v)\le
\phi(u)+\frac{k'}{n+k'}\phi(v).
\]
Applying the decreasing function $\phi^{-1}$ and then
taking the minimum with $v$ gives
\[
C_{n,n+k}(u,v)\ge C_{n,n+k'}(u,v),
\qquad u,v\in[0,1].
\]
Therefore,
\[
(M_n,M_{n+k})\succeq_c(M_n,M_{n+k'}),
\]
as claimed. For the second assertion of the statement, it is enough to observe that the ratio $k/(n+k)$ decreases as the $n$ increases, then reasoning as above.
\end{proof}

\begin{figure}
    \centering
    \includegraphics[width=0.5\linewidth]{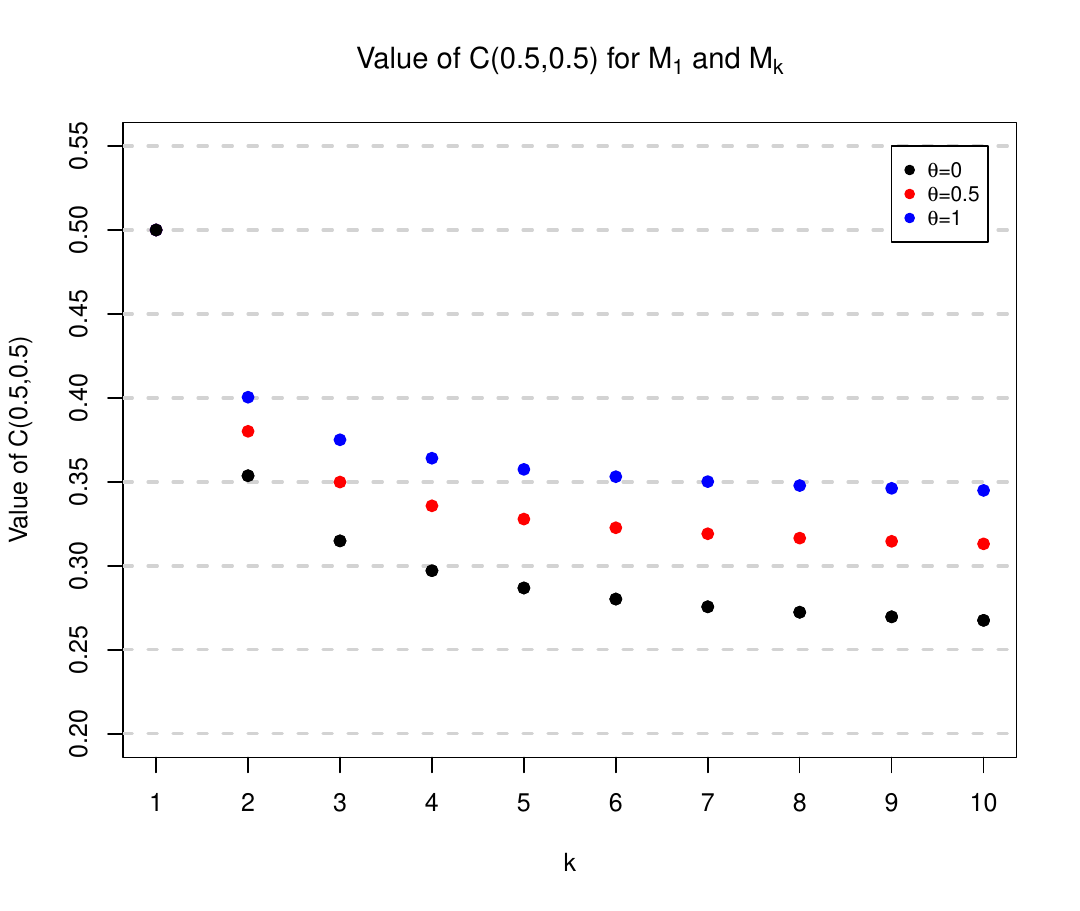}
    \caption{Simulated values of the copula of $(M_1,M_{k})$ in $(0.5,0.5)$ for a sequence of random variables with uniform distribution on $[0,1]$ and a Clayton copula with parameter $\theta\in \{0,0.5,1\}$, using a sample size of $10^6$.}
    \label{fig:3}
\end{figure}

Interestingly, the copula $C_{n,n+k}(u,v)=\min\left\{\phi^{-1}\left(\phi(u)+\phi(v)k/(n+k)\right),v\right\}$ belongs to the class of Archimax copulas~\citep{caperaa2000bivariate} when the function $A$ is defined as $A(t)=\max\left\{t+(1-t)k/(n+k),1-t\right\}$.

The result is illustrated for the Clayton family in  Figure~\ref{fig:3}. However, the result is not true for arbitrary random variables with the same distribution. We provide a counterexample in this regard.

\begin{example}\label{example:fig}
Let $U$ be a standard uniform random variable and let $X_1=U$,
\begin{equation*}
X_2=\begin{cases}
   U+\frac{2}{3} & U\leq \frac{1}{3},\\
   U-\frac{1}{3} & U> \frac{1}{3},\\
   \end{cases} \ \ \textrm{and} \ \ X_3=\begin{cases}
   U+\frac{1}{3} & U\leq \frac{2}{3},\\
   U-\frac{2}{3} & U> \frac{2}{3},\\
\end{cases}
\end{equation*}
see left graph in Figure~\ref{fig:ex}.

Notice that $X_1$, $X_2$ and $X_3$ have standard uniform distribution. Now compute the value of $M_2$ and $M_3$ in terms of the value of $U$,
\begin{equation*}
M_2=\begin{cases}
   U+\frac{2}{3} & U\leq \frac{1}{3},\\
   U & U> \frac{1}{3},\\
\end{cases},\ M_3=\begin{cases}
   U+\frac{2}{3} & U\leq \frac{1}{3},\\
   U+\frac{1}{3} & \frac{1}{3}< U\leq \frac{2}{3},\\
   U & U> \frac{2}{3},\\
\end{cases}
\end{equation*}
see center and right graph in Figure~\ref{fig:ex}. The associated distribution functions are,
\begin{equation*}
    F_2(t)=\begin{cases}
        t-\frac{1}{3} &  \frac{1}{3}\leq t\leq \frac{2}{3},\\
        2(t-\frac{2}{3})+\frac{1}{3} & \frac{2}{3}\leq t\leq 1
    \end{cases},\ F_3(t)=
        3\left(t-\frac{2}{3}\right),\  \frac{2}{3}\leq t\leq 1.
\end{equation*}

Therefore, the $1/3$-quantile of $M_2$ is $2/3$, while $1/3$-quantile of $M_3$ is $\frac{7}{9}$. Compute now some values for the copulas:

\begin{align*}
C_{1,2}\left(\frac{1}{3},\frac{1}{3}\right)=&P\left(X_1\leq \frac{1}{3},M_2\leq F_2^{-1}\left(\frac{1}{3}\right)\right)=P\left(X_1\leq \frac{1}{3},M_2\leq \frac{2}{3}\right)\\=&P\left(U\leq \frac{1}{3},U\in \left(\frac{1}{3},\frac{2}{3}\right]\right)=0.
\end{align*}

\begin{align*}
C_{1,3}\left(\frac{1}{3},\frac{1}{3}\right)=&P\left(X_1\leq \frac{1}{3},M_3\leq F_3^{-1}\left(\frac{1}{3}\right)\right)= P\left(X_1\leq \frac{1}{3},M_3\leq \frac{7}{9}\right)\\=&P\left(U\leq \frac{1}{3},U\in \left[0,\frac{1}{9}\right]\cup \left(\frac{1}{3},\frac{4}{9}\right]\cup \left(\frac{2}{3},\frac{7}{9}\right] \right)P\left(U\leq \frac{1}{9}\right)=\frac{1}{9}.
\end{align*}
so we conclude that $C_{1,2}\not\geq C_{1,3}$, i.e. $(M_1,M_2)\not\succeq_{c} (M_1,M_3)$.

\begin{figure}
    \centering
\includegraphics[width=0.33\textwidth]{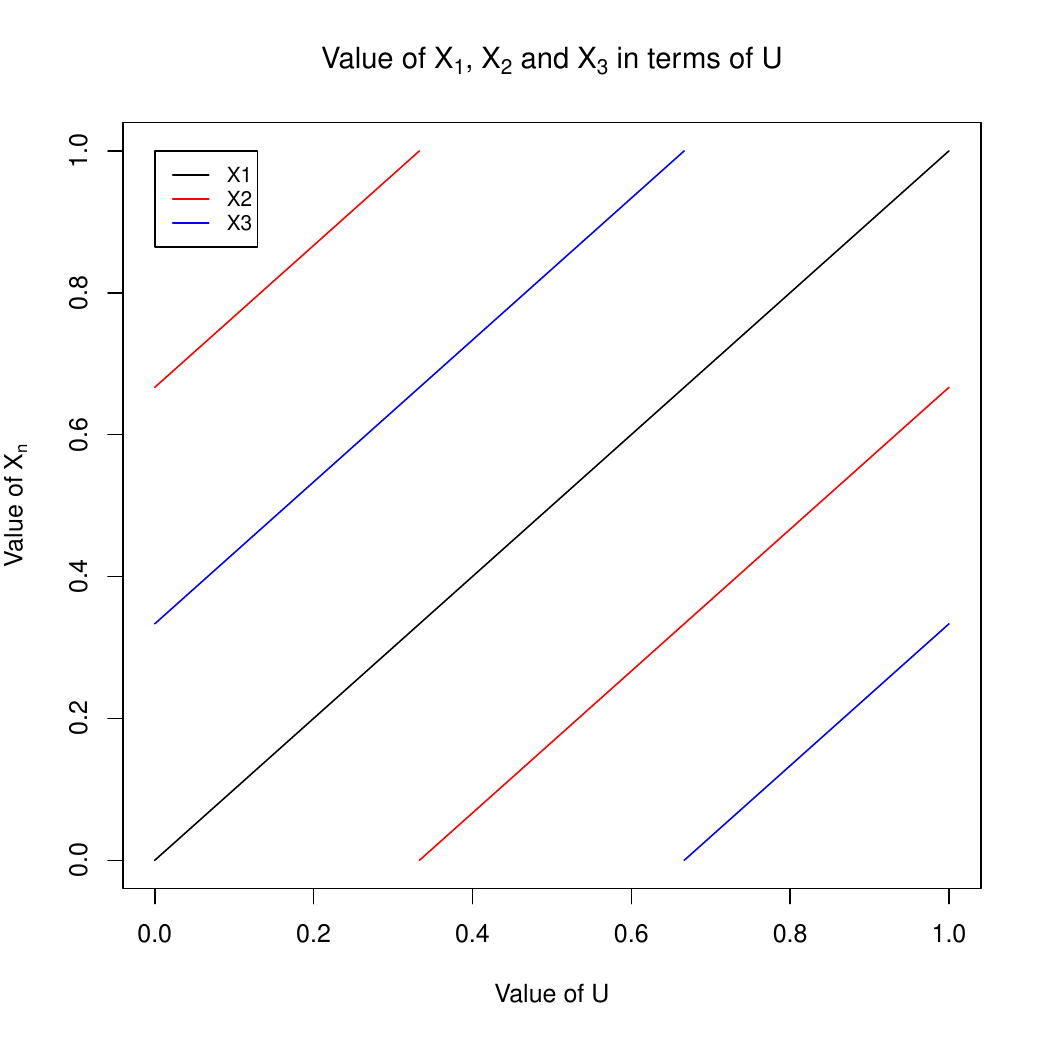}\includegraphics[width=0.33\textwidth]{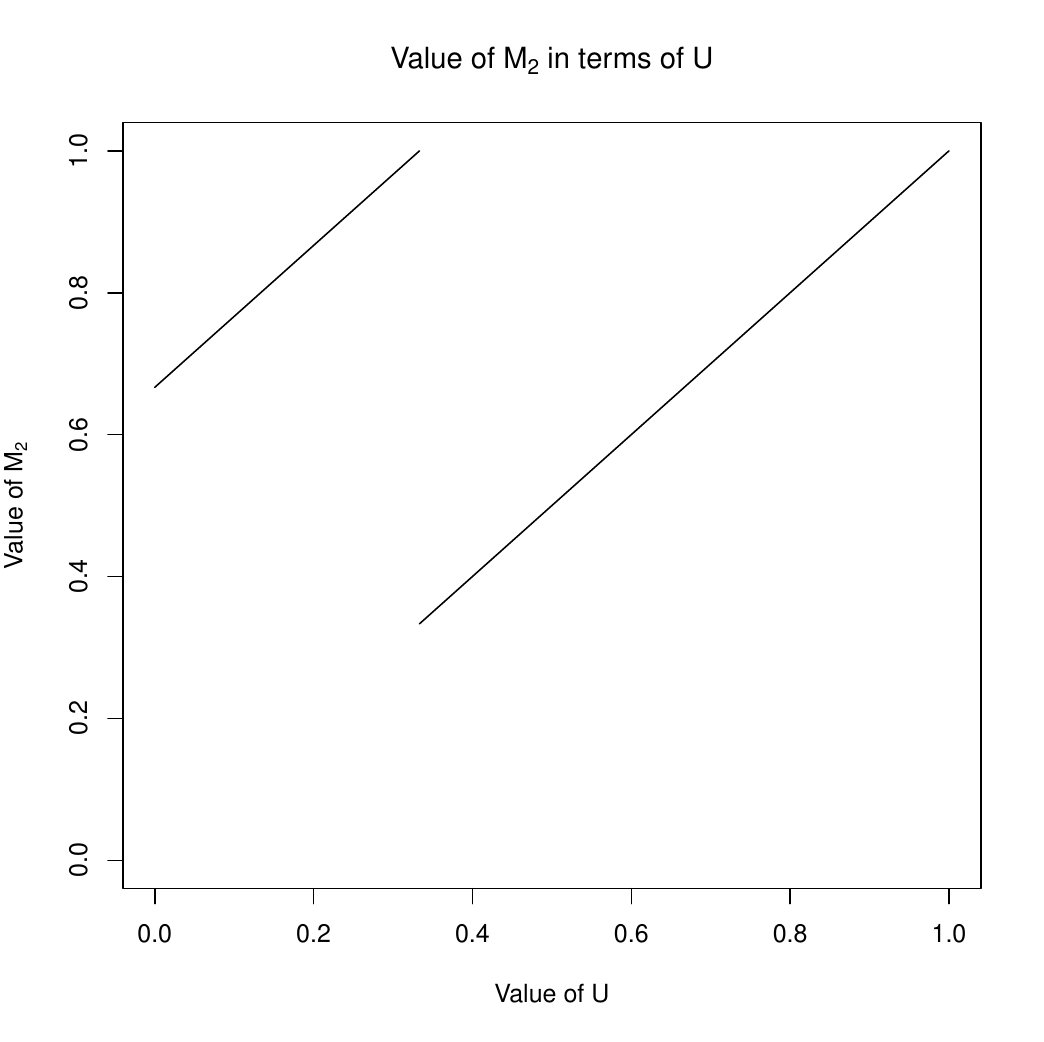}\includegraphics[width=0.33\textwidth]{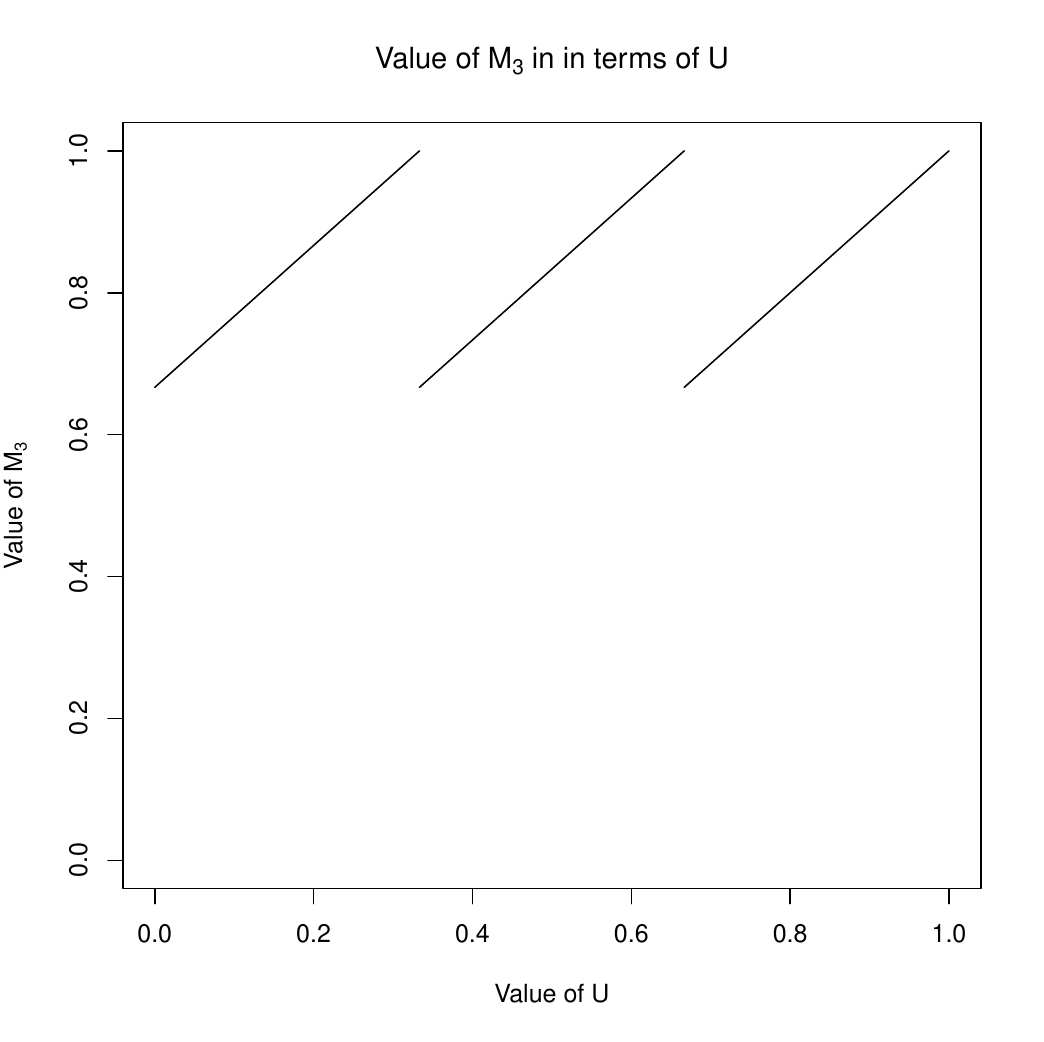}
\caption{Values of $X_1$, $X_2$ and $X_3$ (left), $M_2$ (center) and $M_3$ (right) in terms of the value of $U$ in Example~\ref{example:fig}.}
\label{fig:ex}
\end{figure}

\end{example}

Moving back to Archimedean copulas, for sequences of identically distributed random variables, it is possible to obtain simple formulas for some conditional probabilities. Notice that the distribution function of $M_n$ is $\delta_n(F(x))=\phi^{-1}(n\phi(F(x)))$, where $\phi$ is the generator of the Archimedean copula. For any $x,y\in\mathbb{R}$ such that $x<y$, using the formula of Proposition~\ref{thm:samemarginals_max} we have
\begin{align*}
    P(M_{n+k}\leq y\mid M_{n}\leq x)&=\frac{C_{n,n+k}(F_{n}(x),F_{n+k}(y))}{F_n(x)}\\
    &=\frac{\phi^{-1}\left(\phi(\phi^{-1}(n\phi(F(x))))+\frac{k}{n+k}\phi(\phi^{-1}((n+k)\phi(F(y))))\right)}{\phi^{-1}(n\phi(F(y)))}\\
    &=\frac{\phi^{-1}\left((n\phi(F(x)))+k\phi(F(y))\right)}{\phi^{-1}((n+k)\phi(F(y)))}.
\end{align*}
and
\begin{equation*}
    P(M_{n}\leq x \mid M_{n+k}\leq y)=\frac{C_{n,n+k}(F_{n}(x),F_{n+k}(y))}{F_{n}(y)}=\frac{\phi^{-1}\left((n\phi(F(x)))+k\phi(F(y))\right)}{\phi^{-1}((n+k)\phi(F(y)))}.
\end{equation*}

A natural question that arises is if differences in the dependence structure imply the concordance order between pairs of variables in two different maximum processes (see Figure~\ref{fig:3}). The following result holds directly from Remark~\ref{remark}.

\begin{proposition}\label{prop:conc_dependence}
Let $(X_i,i\in\mathbb{N})$ and $(Y_i,i\in\mathbb{N})$ be two
sequences, each consisting of identically distributed continuous
random variables. Define
\[
M_m=\max(X_1,\ldots,X_m),\qquad\widehat M_m=\max(Y_1,\ldots,Y_m), \qquad m\in\mathbb{N}.
\]
Fix $n,k\in\mathbb{N}$, and let $C$ and $\widehat C$ denote
the copulas of $(X_1,\ldots,X_{n+k})$ and
$(Y_1,\ldots,Y_{n+k})$, respectively. For $i=1,\ldots,n+k$,
define
\[
\delta_i(t)=C(\underbrace{t,\ldots,t}_{i},\underbrace{1,\ldots,1}_{n+k-i}),\qquad\widehat\delta_i(t)=\widehat C(\underbrace{t,\ldots,t}_{i},\underbrace{1,\ldots,1}_{n+k-i}).
\]
For $u,v\in(0,1)$, define
\[
Q(u,v)=C(\underbrace{\delta_n^{-1}(u),\ldots,\delta_n^{-1}(u)}_{n},\underbrace{\delta_{n+k}^{-1}(v),\ldots,\delta_{n+k}^{-1}(v)}_{k})
\]
and
\[\widehat Q(u,v)=\widehat C(\underbrace{\widehat\delta_n^{-1}(u),\ldots,\widehat\delta_n^{-1}(u)}_{n},
\underbrace{\widehat\delta_{n+k}^{-1}(v),\ldots,\widehat\delta_{n+k}^{-1}(v)}_{k}),
\]
where the inverses are generalized inverses.
If
\begin{equation}\label{ineqQ}
    Q(u,v)\le\widehat Q(u,v),
\qquad u,v\in(0,1),
\end{equation}
then
\[
(M_n,M_{n+k})\preceq_c(\widehat M_n,\widehat M_{n+k}).
\]
\end{proposition}

The condition (\ref{ineqQ}) in previous statement can be verified in some cases analytically, or numerically. Moreover, for Archimedean copulas, the representation in Proposition \ref{prop:archi} yields the following sufficient condition in terms of the generators.

\begin{corollary}\label{cor:archi}
Let $(X_i,i\in\mathbb{N})$ and $(Y_i,i\in\mathbb{N})$ be two
sequences, each consisting of identically distributed continuous random variables. Assume that, for every $m\in\mathbb{N}$, the vectors $(X_1,\ldots,X_m)$ and $(Y_1,\ldots,Y_m)$ have Archimedean copulas with generators $\phi$ and $\widehat\phi$, respectively, where each generator is the same in every dimension.

Define
\[
M_m=\max(X_1,\ldots,X_m),
\qquad \widehat M_m=\max(Y_1,\ldots,Y_m),
\qquad m\in\mathbb{N}.
\]
If the function
\[
h(t)=\widehat\phi\bigl(\phi^{-1}(t)\bigr), \qquad t\ge0,
\]
is convex, then
\[
(M_n,M_{n+k})\preceq_c(\widehat M_n,\widehat M_{n+k})\qquad\text{for every }n,k\ge 1.
\]
\end{corollary}

\begin{proof}
Fix $n,k\ge1$ and write
\[
r=\frac{k}{n+k}\in(0,1), \qquad \psi=\phi^{-1},
\qquad \widehat\psi=\widehat\phi^{-1}.
\]
Let $C_{n,k}$ and $\widehat C_{n,k}$ denote the copulas of the two pairs of maxima. By Proposition~\ref{prop:archi},
\[
C_{n,k}(u,v)=\min\left\{v,\, \psi\bigl(\phi(u)+r\phi(v)\bigr)\right\}.
\]
Since $v=\psi(\phi(v))$ and $\psi$ is decreasing, this formula can be rewritten as
\begin{equation}
C_{n,k}(u,v)=\psi\left(\max\{\phi(v),\,\phi(u)+r\phi(v)\}\right).
\label{eq:D}
\end{equation}
Similarly,
\begin{equation}
\widehat C_{n,k}(u,v)=\widehat\psi\left(\max\{\widehat\phi(v),\, \widehat\phi(u)+r\widehat\phi(v)\}\right).
\label{eq:Dhat}
\end{equation}

First consider $u,v\in(0,1)$ and let
\[
a=\phi(u),\qquad b=\phi(v),\qquad \tilde{m}=\max\{b,a+rb\}.
\]
In particular, $\tilde{m}>0$ and
\begin{equation}
0\le a\le \tilde{m},\qquad 0\le b\le \tilde{m}, \qquad a+rb\le \tilde{m}.
\label{eq:bounds}
\end{equation}
The function $h=\widehat\phi\circ\psi$ is increasing and satisfies
$h(0)=\widehat\phi(1)=0$. By convexity, for every $x\in[0,\tilde{m}]$,
\begin{equation}
h(x)=h\left(\frac{x}{\tilde{m}}\tilde{m}+\left(1-\frac{x}{\tilde{m}}\right)0\right)\le \frac{x}{\tilde{m}}h(\tilde{m}).\label{eq:chord}
\end{equation}
Applying \eqref{eq:chord} to $a$ and $b$ and using
\eqref{eq:bounds}, we obtain
\begin{equation*}
h(a)+r h(b)\le \frac{a+rb}{\tilde{m}}h(\tilde{m})\le h(\tilde{m}).
\end{equation*}
Moreover, $h(b)\le h(\tilde{m})$ because $b\le \tilde{m}$ and $h$ is increasing. Therefore,
\begin{equation}
\max\{h(b),\,h(a)+r h(b)\}\le h(\tilde{m}).
\label{eq:key}
\end{equation}

Now observe that $\widehat\phi(u)=h(a)$ and $\widehat\phi(v)=h(b)$. By \eqref{eq:Dhat}, \eqref{eq:key}, and the fact that $\widehat\psi$ is decreasing,
\begin{align*}
\widehat C_{n,k}(u,v)
&=\widehat\psi\left(
\max\{h(b),\,h(a)+r h(b)\}\right)\ge \widehat\psi(h(\tilde{m}))=\widehat\psi\bigl(\widehat\phi(\psi(\tilde{m}))\bigr)=\psi(\tilde{m})=C_{n,k}(u,v),
\end{align*}
where the last equality follows from \eqref{eq:D}. By continuity of copulas, this inequality extends to the boundary of $[0,1]^2$, proving the desired concordance ordering.
\end{proof}

For instance, for the Clayton copula it is easy to check that the concordance order between the copulas of the initial random sequences implies the concordance order between pairs of maxima.

\begin{example}\label{ex:clayton}
Consider two sequences, each consisting of identically distributed continuous random variables, whose finite-dimensional copulas belong to the Clayton family with parameters $\theta$ and $\widehat\theta$, respectively, where
\[
0\le\theta\le\widehat\theta.
\]
Let
\[
M_m=\max(X_1,\ldots,X_m),\qquad \widehat M_m=\max(Y_1,\ldots,Y_m).
\]
Then, for every $n,k\in\mathbb{N}$,
$
(M_n,M_{n+k})\preceq_c(\widehat M_n,\widehat M_{n+k}).
$

To prove it, note that for
\[
\phi_\theta(u)=\frac{u^{-\theta}-1}{\theta}, \qquad 0<\theta\le\widehat\theta,
\]
we have
\[
h(t)
=\phi_{\widehat\theta}\bigl(\phi_\theta^{-1}(t)\bigr)=\frac{(1+\theta t)^{\widehat\theta/\theta}-1}{\widehat\theta}.
\]
Consequently,
\[
h''(t)=(\widehat\theta-\theta)(1+\theta t)^{\widehat\theta/\theta-2}\ge 0.
\]
If $\widehat\theta>\theta$, this function is convex. Thus, under this condition the assumption of Corollary~\ref{cor:archi} is satisfied.
\end{example}

A similar example can be constructed for other families of generators, such as the Gumbel-Hougaard family. However, the concordance order between the copula of the initial sequences does not implies, in general, the concordance order between pairs of variables in the associated maximum processes. We give a counterexample in this regard.

\begin{example}
Let $X_1$, $X_2$ and $X_3$ be three continuous, independent and identically distributed continuous random variables. By Corollary~\ref{cor:iid_conc}, we have that $(M_1,M_2)\preceq_c(M_2,M_3)$. In addition, since they do not share the same copula, $(M_1,M_2)\prec_c(M_2,M_3)$. Now consider three random variables $Y_1$, $Y_2$ and $Y_3$ with the same distribution as $X_1$ such that $Y_1$ and $Y_2$ are comonotone and $Y_3$ is independent of $Y_1$ and $Y_2$. Trivially $(X_1,X_2,X_3)\preceq_c(Y_1,Y_2,Y_3)$, since the copulas and survival copulas are $C_X(u,v,w)=\bar{C}_X(u,v,w)=uvw$ and $C_Y(u,v,w)=\bar{C}_Y(u,v,w)=\min(u,v)w$. However, it is clear that $(\max(Y_1,Y_2),\max(Y_1,Y_2,Y_3))=_{a.s.}(Y_1,\max(Y_1,Y_3))=_{st}(X_1,\max(X_1,X_2))=_{st}(M_1,M_2) \prec_c (M_2,M_3)$ so $(\max(Y_1,Y_2),\max(Y_1,Y_2,Y_3))\prec_c(\max(X_1,X_2),\max(X_1,X_2,X_3))$.
\end{example}

An example of application in predictive value of unusually high peak loads of the statements described in this section is provided here.

\begin{example}\label{ex_nuovo}
Consider a system operating over successive cycles, and let $X_i$ denote the peak load recorded during cycle $i$. Denote its maximum process as
$
M_n=\max(X_1,\ldots,X_n).
$

The loads may share an unobserved operating condition that persists throughout the observation period. For example, this condition may represent the severity of the environment in which the system operates. An explicit model incorporating this feature can be constructed as follows. Let
\[
Z\sim\operatorname{Gamma}(1/\theta,1/\theta), \qquad \theta>0,
\]
where the parameters are the shape and rate, respectively. Conditional on $Z=z$, assume that the loads are independent and identically distributed with
\[
P(X_i\le x\mid Z=z)=
\exp\left[-z\left(\frac{s}{x}\right)^\alpha\right], \qquad x>0,
\]
for some $s,\alpha>0$. A larger value of $Z$ produces stochastically larger loads. Since
$
E[e^{-tZ}]=(1+\theta t)^{-1/\theta},
$
the common marginal distribution is
\[
F(x)=\left[1+\theta\left(\frac{s}{x}\right)^\alpha\right]^{-1/\theta},\qquad x>0.
\]

Moreover,
\[
P(X_1\le x_1,\ldots,X_m\le x_m)=\left[1+\theta\sum_{i=1}^{m}\left(\frac{s}{x_i}\right)^\alpha \right]^{-1/\theta}.
\]

Thus, for every $m$, the vector $(X_1,\ldots,X_m)$ has a Clayton copula with parameter $\theta$. In particular, Proposition~\ref{prop:archi} applies to this sequence.

Suppose that, after $n$ cycles, the maximum process is unusually high. It is then natural to ask how informative is this observation about the running maximum after another $k$ cycles. On this aim, let $G_m$ denote the distribution function of $M_m$, and, for
$q\in(0,1)$, define
\[
p_{n,k}(q)=P\left(M_{n+k}>G_{n+k}^{-1}(q)\;\middle|\;M_n>G_n^{-1}(q)\right).
\]

Here, an unusually high maximum means one exceeding the $q$-quantile at the corresponding horizon. Both exceedance events therefore have unconditional probability $1-q$.

If $C_{n,n+k}$ denotes the copula of $(M_n,M_{n+k})$, then
\[
p_{n,k}(q)=\frac{1-2q+C_{n,n+k}(q,q)}{1-q}.
\]

Consequently, Proposition~\ref{prop:archi} yields
\[
p_{n,k}(q)\ge p_{n,k'}(q)\qquad\text{whenever }k\le k'.
\]

An unusually high current maximum is therefore a stronger predictor of an unusually high future maximum over shorter horizons. For example for the Clayton model, writing $r_{n,k}=k/(n+k)$, we have
\[
C_{n,n+k}(u,v)=\min\left\{v,\ \left[u^{-\theta}+r_{n,k}(v^{-\theta}-1)\right]^{-1/\theta}\right\}.
\]

Hence,
\[
p_{n,k}(q)=\frac{1-2q+\left[q^{-\theta}+r_{n,k}(q^{-\theta}-1)\right]^{-1/\theta}}{1-q}.
\]

For illustration, take $\theta=1$, $n=30$, and $q=0.95$. In this case, 
\[
p_{n,k}(q)=
\frac{1+r_{n,k}(1-2q)}{1+r_{n,k}(1-q)},
\]
giving the following probabilities:
\begin{gather*}
p_{30,1}(0.95)=0.969,\quad p_{30,5}(0.95)=0.865,\quad p_{30,10}(0.95)=0.765,\\ p_{30,30}(0.95)=0.537,\quad p_{30,90}(0.95)=0.313,\quad p_{30,\infty}(0.95)=0.095.
\end{gather*}
where $p_{30,\infty}(0.95)$ denotes the limiting probability.

Thus, conditional on the current maximum exceeding its own $95$th percentile, the probability that the future maximum exceeds its corresponding $95$th percentile decreases with the forecast horizon. Nevertheless, its limit remains above the unconditional probability of $0.05$.

This behavior reflects two sources of dependence. At short horizons, an abnormally big value of the maximum process at a particular time implies, since the process is almost surely increasing, abnormally big values for the subsequent times. At long horizons, the common operating condition $Z$ continues to provide predictive information, since the variables in the initial sequence $(X_n,n\in\mathbb{N})$ have a dependence that does not decay with time.
\end{example}

We want to point out that the latter discussion is not specific to the chosen marginal distribution or to the Clayton family. It holds for any sequence satisfying the assumptions of Proposition~\ref{prop:archi}.

\subsection{Asymptotic results}

\begin{figure}
    \centering
    \includegraphics[width=0.45\linewidth]{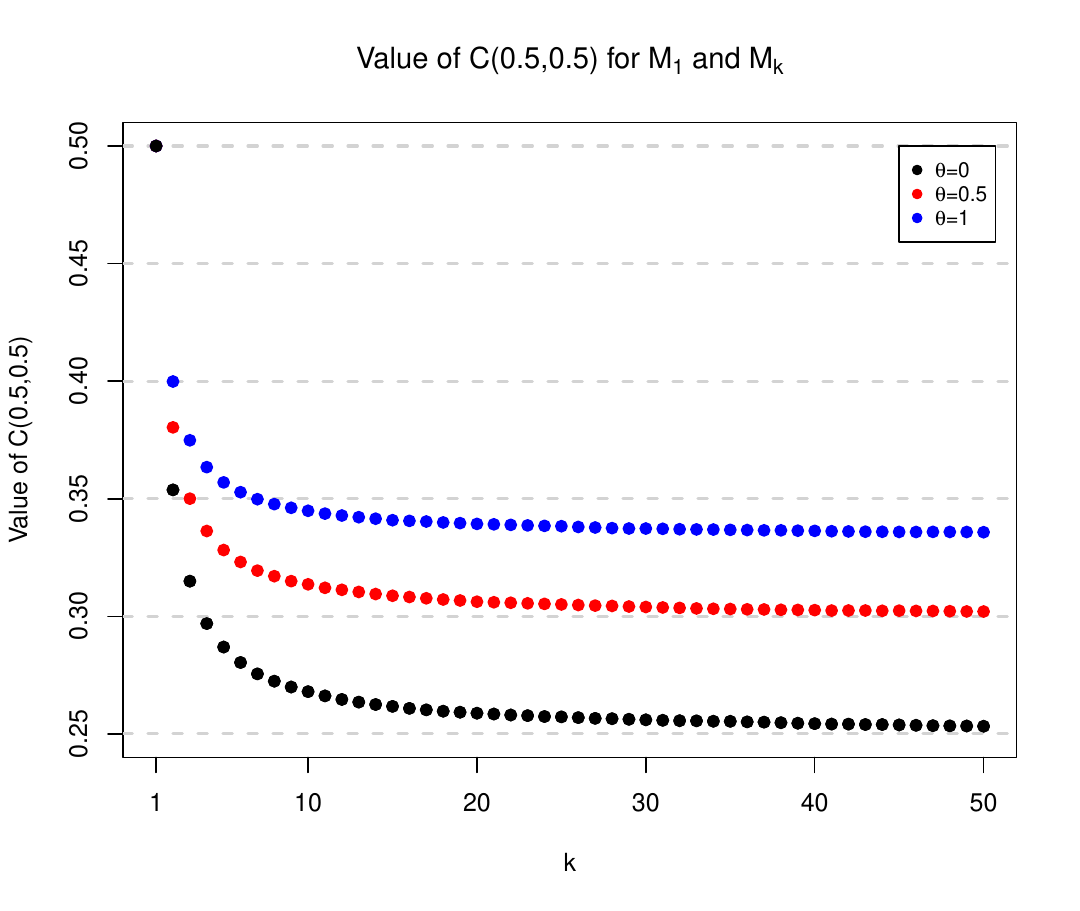}
    \includegraphics[width=0.45\linewidth]{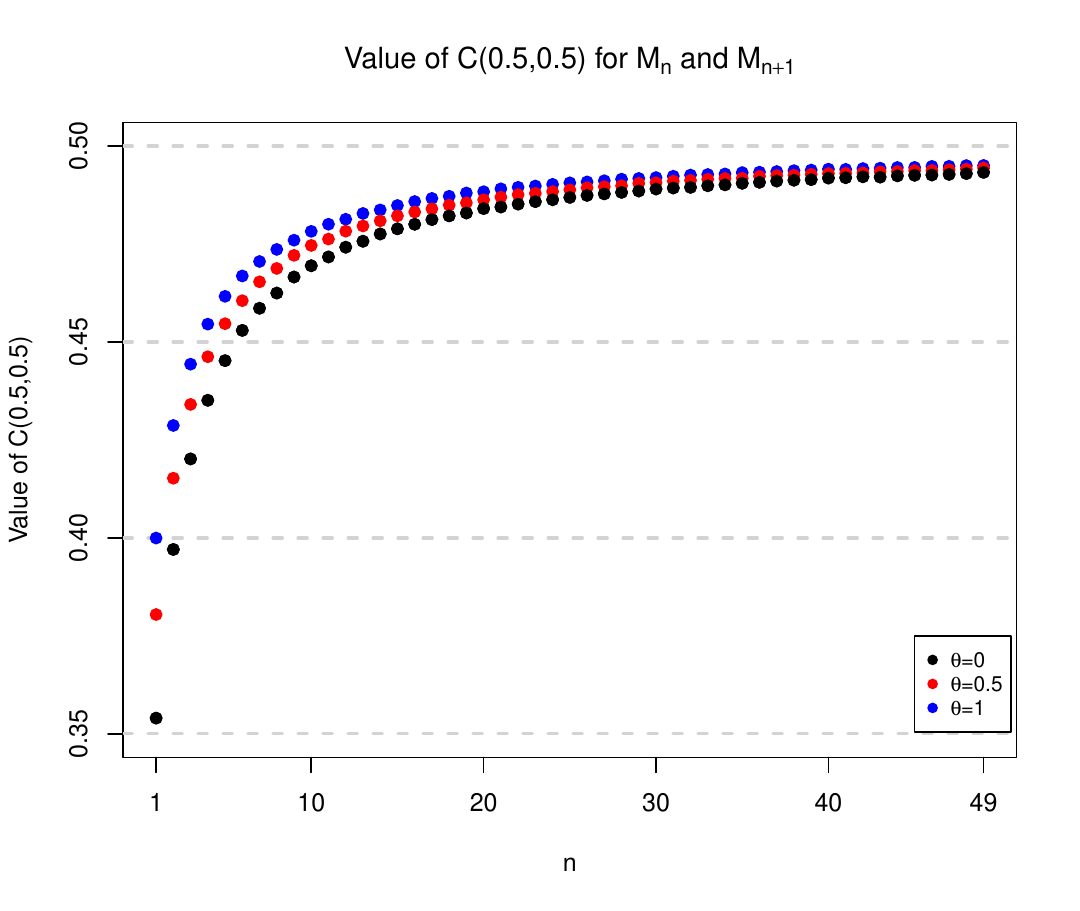}
    \caption{Simulated values of the copula of $(M_1,M_{k})$ (left) and $(M_{n},M_{n+1})$ in $(0.5,0.5)$ for a sequence of random variables with uniform distribution on $[0,1]$ and a Clayton copula with parameter $\theta\in \{0,0.5,1\}$, using a sample size of $10^6$.}
    \label{fig:4}
\end{figure}

Keeping in mind Figures~\ref{fig:1},~\ref{fig:2} and~\ref{fig:3}, it seems that there could be a convergence of the copula of $M_n$ and $M_{n+k}$ as $n$ or $k$ tend to infinity. Simulating bigger values for them, see Figure~\ref{fig:4}, that tendency is more clear.

This section is devoted to follow that intuition and study the asymptotic behaviors of the latter copulas. It is worth to note that the asymptotic inequalities already proved imply the convergence of the associated copulas under some conditions. If $(C_n,n\in\mathbb{N})$ is a sequence of copulas such that $C_n(u,v)\leq[\geq] C_{n+1}(u,v)$ for any $u,v\in [0,1]$ and $n\in \mathbb{N}$, then the sequence $(C_n(u,v),n\in\mathbb{N})$ is increasing [decreasing] and bounded, thus convergent. However, the value of the limit remains unknown. We will focus on this task in the rest of the section.

Starting from the case of independent variables, we expect the random variables $M_n$ and $M_{n+k}$ to be independent when $k$ tends to infinity if the sequence increases their values with time.

\begin{proposition}\label{prop:conv_ind_1}
Let $(X_n,n\in\mathbb{N})$ be a sequence of continuous independent random variables such that $X_i\leq_{st} X_j$ for any $i,j\in\mathbb{N}$ with $i\leq j$. Then,
\begin{equation*}
    \lim_{k\to\infty}C_{n,n+k}(u,v)=uv
\end{equation*}
for any $n\in\mathbb{N}$.
\end{proposition}

\begin{proof}
Fix $n\ge1$ and $v\in(0,1)$, and define $q_k=G_{n+k}^{-1}(v)$. Since $G_{n+k}$ is continuous,
$
G_{n+k}(q_k)=v.
$

The assumed stochastic ordering implies
\[
F_j\le F_n \quad\text{for }j>n,\qquad F_i\ge F_n \quad\text{for }i\le n.
\]
By independence,
\[
v=\prod_{i=1}^{n+k}F_i(q_k)=\left(\prod_{i=1}^{n}F_i(q_k)\right)\left(\prod_{i=n+1}^{n+k}F_i(q_k)\right)\leq \prod_{i=n+1}^{n+k}F_i(q_k)\leq \prod_{i=n+1}^{n+k}F_n(q_k)=F_n(q_k)^k.
\]
Consequently,
\[
1\ge G_n(q_k)=\prod_{i=1}^{n}F_i(q_k)\ge F_n(q_k)^n = \left(F_n(q_k)^k\right)^{n/k}=v^{n/k}.
\]
It follows that
\[
\lim_{k\to\infty} G_n\bigl(G_{n+k}^{-1}(v)\bigr)=1
\]

Using Proposition~\ref{thm:copula_ind}, for $u,v\in(0,1)$ we obtain
\[
\lim_{k\to\infty} C_{n,n+k}(u,v)=\lim_{k\to\infty} \min\left\{\frac{uv}{G_n(G_{n+k}^{-1}(v))},v\right\}=uv,
\]
where in the last step we have used the continuity of the involved functions. The same conclusion holds on the boundary by the boundary
conditions of copulas.
\end{proof}

Notice that the latter result also holds for a sequence of independent and identically distributed random variables.  In this sense, the result is similar to the well-known convergence of the copula of the minimum and the maximum of independent and identically distributed random variables, see~\cite{schmitz2004revealing}.  If we drop independence the result is not longer true and the limit might not exist. A counterexample in this regard is given below; here the continuity and stochastic
monotonicity assumptions of Proposition~\ref{prop:conv_ind_1} are satisfied, and even almost sure monotonicity, but not independence.

\begin{example}\label{ex:no limit}
Let $U$ be uniformly distributed on $[0,1]$, and define
\[
X_n=\begin{cases}
n-1+U, & n \text{ odd},\\
n-U,   & n \text{ even}.
\end{cases}
\]
Then $X_n$ is uniformly distributed on $[n-1,n]$ for every
$n\in\mathbb{N}$. In particular,
\[
X_i\le_{\mathrm{st}}X_j\qquad\text{whenever }i\le j.
\]

Moreover,
\[
X_{n+1}-X_n=\begin{cases}
2(1-U), & n \text{ odd},\\
2U,     & n \text{ even},
\end{cases}
\]
so the sequence is almost surely nondecreasing. Consequently,
\[
M_n=\max(X_1,\ldots,X_n)=_{a.s.}X_n.
\]

If $k$ is even, $X_n$ and $X_{n+k}$ are both increasing functions of $U$ or both decreasing functions of $U$. They are therefore comonotone. If $k$ is odd, one is an increasing function of $U$ and the other is a decreasing function of $U$, so they are countermonotone. Thus, the copula of $(M_n,M_{n+k})$ is
\[
C_{n,n+k}(u,v)=\begin{cases}
\min(u,v),       & k \text{ even},\\
\max(u+v-1,0),   & k \text{ odd}.
\end{cases}
\]
In particular,
\[
C_{n,n+k}\left(\frac12,\frac12\right)=\begin{cases}
\frac12, & k \text{ even},\\
0,       & k \text{ odd}.
\end{cases}
\]
Hence, for every fixed $n$, the sequence of copulas
$(C_{n,n+k},k\in\mathbb{N})$ does not converge as $k$ goes to infinity.
\end{example}

Regarding the limit in $n$, one should expect the probability of having an update should go to 0 and, therefore, to observe a strong positive dependence of $(M_n,M_{n+k})$ for large values of $n$. However, to prove it formally we must add the assumption of the uniform convergence of the sequence of the distribution functions.

\begin{proposition}\label{prop:uniform-conv-min}
Let $(X_n,n\in\mathbb{N})$ be a sequence of independent random variables that converges, in distribution, to a variable $X$ having continuous distribution $F$. Denote the distribution function of $X_i$ as $F_i$ with $i\in\mathbb{N}$ and  the copula of $(M_n,M_{n+k})$ as $C_{n,n+k}$. Then, for every fixed positive integer $k$,
\[
\lim_{n\to\infty}\sup_{(u,v)\in[0,1]^2}\left|C_{n,n+k}(u,v)-\min(u,v)\right|=0,
\]
and therefore also
\begin{equation*}
    \lim_{n\to\infty}C_{n,n+k}(u,v)=\min(u,v).
\end{equation*}
\end{proposition}

\begin{proof}
Fix $k\ge1$, and let $p_{n,k}=P(M_{n+k}>M_n).$ We first show that $p_{n,k}\to 0$.

For every $t\in\mathbb{R}$,
\[
\{M_{n+k}>M_n\}\subseteq\{M_n\le t\}\cup\bigcup_{j=1}^{k}\{X_{n+j}>t\}.
\]
Indeed, if $M_n>t$ and none of the next $k$ observations exceeds
$t$, the running maximum cannot increase. Consequently,
independence gives
\begin{equation}\label{eq:prop38_update_bound}
p_{n,k}\le\prod_{i=1}^{n}F_i(t)+\sum_{j=1}^{k}\bigl[1-F_{n+j}(t)\bigr].
\end{equation}

Since $F$ is continuous, convergence in distribution implies $F_i(t)\longrightarrow F(t)$ for every $t\in\mathbb{R}$.

Fix $t$ such that $F(t)<1$, and choose $\rho\in(F(t),1)$.
Since there exists $i_0$ such that $F_i(t)\le\rho$ for every
$i\ge i_0$ then, for $n\ge i_0$,
\[
0\le\lim_{n\to\infty} \prod_{i=1}^{n}F_i(t)\le\lim_{n\to\infty} \rho^{\,n-i_0+1}=0.
\]
Taking the limit superior in \eqref{eq:prop38_update_bound}, we obtain
\[
\limsup_{n\to\infty}p_{n,k}\le k[1-F(t)].
\]

Since the limit of distribution functions to a continuous distribution function is always a uniform limit, see Lemma 2.11 in~\cite{van2000asymptotic}, choose $t$ such that $F(t)<1$ is arbitrarily close to one. Therefore,
\begin{equation}\label{eq:prop38_update_vanishes}
\lim_{n\to\infty} p_{n,k}= 0.
\end{equation}

We next relate this probability to the distance between
$C_{n,n+k}$ and the comonotonic copula. Let
\[
G_m(x)=P(M_m\le x)=\prod_{i=1}^{m}F_i(x).
\]
Each $G_m$ is continuous. Fix $u,v\in(0,1)$, and set
\[
a=G_n^{-1}(u),\qquad b=G_{n+k}^{-1}(v),
\]
where the inverses are generalized inverses. By continuity,
\[
G_n(a)=u,\qquad G_{n+k}(b)=v.
\]

On the event $\{M_n=M_{n+k}\}$, the inequality $M_n\le\min(a,b)$ implies both $M_n\le a$ and $M_{n+k}\le b$. Thus,
\begin{align*}
C_{n,n+k}(u,v)
&=P(M_n\le a,\ M_{n+k}\le b)\ge
P\bigl(M_n\le\min(a,b),\ M_n=M_{n+k}\bigr)\\
&\ge G_n(\min(a,b))-p_{n,k}=\min\{u,G_n(b)\}-p_{n,k}.
\end{align*}
Since $M_n\le M_{n+k}$ almost surely,
\[
G_n(b)\ge G_{n+k}(b)=v.
\]
It follows that
\[
C_{n,n+k}(u,v)\ge\min(u,v)-p_{n,k}.
\]
Combining this inequality with the Fr\'echet--Hoeffding upper bound and extending to the boundary by continuity yields
\begin{equation}
0\le\min(u,v)-C_{n,n+k}(u,v) \le p_{n,k},
\qquad (u,v)\in[0,1]^2.
\end{equation}
The desired uniform convergence now follows from~\eqref{eq:prop38_update_vanishes}.
\end{proof}

In the following we give an example in which, if the limit is not continuous, the convergence to the comonotone copula is not longer true.

\begin{example}\label{ex:no_uniform} 
Let $(X_n,n\in\mathbb{N})$ be a sequence of independent random variables with distribution functions
\[
F_n(x)=\begin{cases}
0,      & x<0,\\
x^{n!}, & 0\le x\le1,\\
1,      & x>1.
\end{cases}
\]
Note that the sequence $(F_n, n \in \mathbb{N})$ of continuous distributions converges pointwise to the distribution function of a point mass at $1$, but the convergence is not uniform.

Define $S_n=\sum_{i=1}^{n}i!$. By independence, the distribution function of
$M_n=\max(X_1,\ldots,X_n)$ satisfies
\[
G_n(x)=\prod_{i=1}^{n}F_i(x)=x^{S_n},\qquad 0\le x\le1.
\]
Consequently,
\[
G_{n+k}^{-1}(v)=v^{1/S_{n+k}},\qquad\textrm{and}\qquad G_n\bigl(G_{n+k}^{-1}(v)\bigr)=v^{S_n/S_{n+k}}\qquad v\in(0,1).
\]

Fix $k\ge1$. Using $S_n\le2n!$ and
$S_{n+k}\ge(n+k)!$, we obtain
\[
0\le\lim_{n\to \infty}\frac{S_n}{S_{n+k}}\le\lim_{n\to \infty}\frac{2n!}{(n+k)!}\le\lim_{n\to \infty}\frac{2}{n+1}=0.
\]
Therefore, for every $v\in(0,1)$,
\[
\lim_{n\to \infty} G_n\bigl(G_{n+k}^{-1}(v)\bigr)=1.
\]

By Proposition~\ref{thm:copula_ind}, the copula of $(M_n,M_{n+k})$ is
\[
C_{n,n+k}(u,v)=\min\left\{u\,v^{\,1-S_n/S_{n+k}},\,v\right\},\qquad u,v\in(0,1).
\]
Hence,
\[
\lim_{n\to\infty}C_{n,n+k}(u,v)=\min\{uv,v\}=uv.
\]
The same limit holds on the boundary by the boundary conditions of copulas. Thus, for every fixed $k\ge1$, the limiting copula is the independence copula rather than the comonotone copula. In particular, convergence of the individual distribution functions does not suffice to ensure the conclusion of Proposition~\ref{prop:uniform-conv-min}.
\end{example}

In the case of exchangeable sequences, it is also possible to prove that, under appropriate conditions, the limit copula when $n$ tends to infinity is the minimum copula. This is shown in the following statement.

\begin{proposition}\label{prop:finetti-min}
Let $(X_i,i\in\mathbb{N})$ be an exchangeable sequence of continuous random variables, and define
\[
M_n=\max(X_1,\ldots,X_n).
\]
Let $C_{n,n+k}$ denote the copula of $(M_n,M_{n+k})$. Then, for every $n,k\in\mathbb{N}$,
\[
0\le\min(u,v)-C_{n,n+k}(u,v)\le\frac{k}{n+k},\qquad u,v\in[0,1].
\]
Consequently, for every fixed $k\in\mathbb{N}$,
\[\lim_{n\to\infty}\sup_{u,v\in[0,1]}\left|C_{n,n+k}(u,v)-\min(u,v)\right|=0
\]
and therefore also
\begin{equation*}
    \lim_{n\to\infty}C_{n,n+k}(u,v)=\min(u,v).
\end{equation*}

More generally, the uniform convergence holds with
$k=k_n$ whenever $(k_n)$ is a sequence of positive integers
such that $\lim_{n\to\infty}k_n/n=0$.
\end{proposition}

\begin{proof}
We first bound the probability that the running maximum increases between times $n$ and $n+k$. Conditional on $(X_1,\ldots,X_{n+k})$, select an index $J$ uniformly at random from the nonempty set
\[
\left\{i\in\{1,\ldots,N\}:X_i=M_{n+k}\right\}.
\]

By exchangeability and the symmetry of the selection rule, $J$ is uniformly distributed on $\{1,\ldots,n+k\}$. On the event $\{M_{n+k}>M_n\}$, every index attaining $M_{n+k}$ belongs to $\{n+1,\ldots,n+k\}$. Hence,
\[
\{M_{n+k}>M_n\}\subseteq\{J>n\},
\]
and therefore
\begin{equation}\label{eq:maximum_update_bound}
p_{n,k}=P(M_{n+k}>M_n)
\le P(J>n)
=\frac{k}{n+k}.
\end{equation}

Fix $u,v\in(0,1)$ and define $a=G_n^{-1}(u)$ and $b=G_{n+k}^{-1}(v)$, where the inverses are generalized inverses. Continuity gives
\[
G_n(a)=u,\qquad G_{n+k}(b)=v.
\]

On the event $\{M_n=M_{n+k}\}$, the condition
$M_n\le\min(a,b)$ implies both $M_n\le a$ and $M_N\le b$. Consequently,
\begin{align*}
C_{n,n+k}(u,v)
&=P(M_n\le a,M_{n+k}\le b)\ge
P\bigl(M_n\le\min(a,b),\,M_n\\
&=M_{n+k}\bigr)\ge
G_n(\min(a,b))-p_{n,k}=
\min\{u,G_n(b)\}-p_{n,k}.
\end{align*}
Since $M_n\le M_{n+k}$ almost surely, we have $G_n(b)\ge G_{n+k}(b)=v.$ It follows that
\[
C_{n,n+k}(u,v)\ge\min(u,v)-p_{n,k}\ge\min(u,v)-\frac{k}{n+k},
\]
where the last inequality follows from
\eqref{eq:maximum_update_bound}. Combining this lower bound with the upper Fr\'echet--Hoeffding bound,
\[
C_{n,n+k}(u,v)\le\min(u,v),
\]
and extending to the boundary by continuity proves
\[
0\le\min(u,v)-C_{n,n+k}(u,v)\le\frac{k}{n+k}
\]
throughout $[0,1]^2$. The uniform convergence for fixed $k$ follows immediately. If $k=k_n$ and $\lim_{n\to\infty} k_n/n=0$, then
\[
\lim_{n\to \infty} \frac{k_n}{n+k_n}=0,
\]
which proves the final assertion.
\end{proof}

Unfortunately, the limit when $k$ tends to infinite seems unaffordable for general exchangeable sequences. However, if one considers an Archimedean copula, the limits are easy to compute.

\begin{proposition}\label{prop:limit-archi}
Let $(X_n,n\in\mathbb{N})$ be a sequence of identically distributed continuous random variables such that, for any $n\in\mathbb{N}$, $(X_1,\dots,X_n)$ has an Archimedean copula with generator $\phi$. Then, denoting as $C_{n,n+k}$ the copula of $(M_n,M_{n+k})$, it holds that
\begin{equation*}
    \lim_{n\to\infty}C_{n,n+k}(u,v)=\min(u,v),
\end{equation*}
for any $k\in\mathbb{N}$ and
\begin{equation*}
    \lim_{k\to\infty}C_{n,n+k}(u,v)=C(u,v),
\end{equation*}
for any $n\in\mathbb{N}$, with $C$ being the copula of $(X_1,X_2)$.
\end{proposition}

\begin{proof}
Using the formulas derived in the proof of Proposition~\ref{prop:archi}, we have that the expression of the copula is
\begin{equation*}
    C_{n,n+k}(u,v)=\begin{cases}
        v & \text{if } \phi^{-1}\left(\frac{1}{n}\phi(u)\right)\geq \phi^{-1}\left(\frac{1}{n+k}\phi(v)\right),\\
\phi^{-1}\left(\phi(u)+\frac{k}{n+k}\phi(v)\right) & \text{if } \phi^{-1}\left(\frac{1}{n}\phi(u)\right)< \phi^{-1}\left(\frac{1}{n+k}\phi(v)\right)
    \end{cases}
\end{equation*}
which, using that $\phi^{-1}$ is strictly decreasing, is equivalent to
\begin{equation*}
    C_{n,n+k}(u,v)=\begin{cases}
        v & \text{if } \frac{1}{n}\phi(u)\leq \frac{1}{n+k}\phi(v),\\
\phi^{-1}\left(\phi(u)+\frac{k}{n+k}\phi(v)\right) & \text{if } \frac{1}{n}\phi(u)> \frac{1}{n+k}\phi(v).
    \end{cases}
\end{equation*}

Then, using the continuity of the associated functions and the fact that $\phi$ is strictly decreasing, we have
\begin{align*}
   \lim_{n\to\infty} C_{n,n+k}(u,v)=&\lim_{n\to\infty} \begin{cases}
        v & \text{if } \frac{1}{n}\phi(u)\leq \frac{1}{n+k}\phi(v)\\
\phi^{-1}\left(\phi(u)+\frac{k}{n+k}\phi(v)\right) & \text{if } \frac{1}{n}\phi(u)> \frac{1}{n+k}\phi(v)
    \end{cases}\\=& \begin{cases}
        v & \text{if } u\geq v\\
\phi^{-1}\left(\phi(u)\right) & \text{if } u<v
    \end{cases}=\begin{cases}
        v & \text{if } u\geq v\\
u & \text{if } u<v
    \end{cases}=\min(u,v),
\end{align*}
and
\begin{align*}
\lim_{k\to\infty} C_{n,n+k}(u,v)=&\lim_{k\to\infty} \begin{cases}
        v & \text{if } \frac{1}{n}\phi(u)\leq \frac{1}{n+k}\phi(v)\\
\phi^{-1}\left(\phi(u)+\frac{k}{n+k}\phi(v)\right) & \text{if } \frac{1}{n}\phi(u)> \frac{1}{n+k}\phi(v)
    \end{cases}\\=& \begin{cases}
        v & \text{if } \frac{1}{n}\phi(u)\leq 0\\
\phi^{-1}\left(\phi(u)+\phi(v)\right) & \text{if } \frac{1}{n}\phi(u)> 0
    \end{cases}=\phi^{-1}\left(\phi(u)+\phi(v)\right),
\end{align*}
which is the copula of $(X_1,X_2)$.
\end{proof}

It is worth to mention that, since Archimedean copulas are symmetric, the sequence $(X_n,n\in\mathbb{N})$ is exchangeable, so $(X_i,X_j)$ have the same copula for any $i,j\in\mathbb{N}$. This result is coherent with Figure~\ref{fig:4}, in particular we have that bivariate Clayton copulas with parameters $\theta=0$, $\theta=0.5$ and $\theta=1$ equal, respectively, $0.25$, $0.299$ and $0.333$ in $(0.5,0.5)$. Moreover, the points in the right plot go closer to $0.5$ (and $\min(0.5,0.5)=0.5$).

\section{Conclusion and open problems}

In this paper, the dependence properties of the discrete-time maximum process $(M_n,n\in\mathbb{N
})$ is studied considering a range of assumptions for the initial random sequence $(X_n,n\in\mathbb{N})$, including independence, same distribution among the variables or exchangeability. Under different conditions it has been proved that the pair $(M_n,M_{n+k})$ is decreasing in the concordance order as $k$ increases and is increasing in the concordance order as $n$ increases. Under appropriate assumptions, the limit case is also studied. If $n$ tends to infinity, then the copula converges to the minimum copula and, if $k$ tends to infinity, the copula converges to the copula of consecutive random variables in $(X_n,n\in\mathbb{N})$.

An interesting question for further research is to study the application of this results to extreme distribution theory. The here-presented results provide limit behavior for the copula, which combined with the already known limits for the scaled marginals~\citep{leadbetter2012extremes}, could give bivariate limit distributions. It should be studied in detail the conditions under which the limit of the bivariate distribution function can be written as the limit of the copula compose with the limit of the marginals and how the dependence in the initial sequence affects the marginal distributions, see~\cite{herrmann2024limiting} for some information in this regard.

\section*{Funding}
J. Baz is supported by the Ministry of Science and Innovation (PID2022-139886NB-I00) and the Principality of Asturias (SEK-25-GRU-GIC-24-018). F. Pellerey is partially supported by the research group GNAMPA of INdAM (Istituto Nazionale di Alta Matematica).

\bibliographystyle{abbrv}

\bibliography{biblio.bib}

\end{document}